\documentclass[10pt]{article}

\usepackage[utf8]{inputenc}
\usepackage[english]{babel}
\usepackage[square]{natbib}
\usepackage{hyperref}
\usepackage{xcolor}
\usepackage{tikz}
\hypersetup{
    colorlinks=false,    
    pdfborderstyle={},   
    pdfborder={0 0 0},  
    linkcolor=black,
    citecolor=black,
    filecolor=black,
    urlcolor=black
}

\usepackage{amsmath, amssymb, bm}

\usepackage{graphicx}
\usepackage{booktabs}
\usepackage{multirow}
\usepackage{makecell}
\usepackage{enumitem}

\usepackage[margin=1.5in]{geometry}
\makeatletter
\renewcommand\footnoterule{%
  \vspace{0.7em}   
  \hrule width 0.4\columnwidth height 0.4pt
  \kern 0.3em      
}
\makeatother
\usepackage[hang]{footmisc}
\usepackage[utf8]{inputenc}
\usepackage[english]{babel}
\usepackage{amsmath}
\usepackage{booktabs}
\usepackage{float}
\usepackage{makecell}
\usepackage{tabularx} 

\usepackage[colorinlistoftodos]{todonotes}

\begin{document}

\definecolor{frontMinMax}{HTML}{2E8B57}
\definecolor{frontMinMin}{HTML}{8C564B}
\definecolor{frontMaxMax}{HTML}{E07B39}
\definecolor{frontMaxMin}{HTML}{7B2CBF}
\title{\textbf{Ordering Is Not Deciding: Towards Preference-Consistent Multi-Objective Design--Decision Optimisation}}

\renewcommand{\thefootnote}{\fnsymbol{footnote}}

\author{
E. \"Ozer$^{1}$, A.R.M. Wolfert$^{2,}\footnote{Corresponding author: A.R.M. (Rogier) Wolfert, \texttt{a.r.m.wolfert@tudelft.nl}}$\\[0.5em]  
\small$^{1}$ Decision Science \& Engineering, Boskalis, The Netherlands\\[-0.3em]  
\small$^{2}$ Department of Algorithmics, Faculty of Mathematics and Computer Science\\[-0.3em]  
\small Delft University of Technology, The Netherlands
}

\renewcommand{\thefootnote}{\arabic{footnote}}
\date{}
\maketitle

{
\renewcommand{\baselinestretch}{0.96}\normalsize
\vspace*{-2em}

\begin{abstract}

\noindent
Contemporary integrated design and group decision-making problems are frequently addressed using Pareto analysis. However, Pareto filtering merely orders the solution space by dominance to identify a set of non-dominated alternatives. Identifying a best-fit design--decision solution from this set requires additional actor-preference information and an \emph{a posteriori} selection step. Pareto analysis therefore cannot solve a design--decision problem without preference information. A further fundamental question then arises: even when preference information is applied \emph{a posteriori}, does Pareto filtering guarantee that the best-fit solution is retained?

\noindent
This paper demonstrates that the answer to this question is no. When preferences are introduced only \emph{a posteriori}, dominance filtering may already have excluded the alternative with the highest aggregated actor preference. Consequently, selecting the most-preferred alternative from the resulting Pareto set cannot guarantee identification of the best-fit solution, even when all combinations of fixed performance-improvement directions are considered. Pareto filtering may therefore exclude the best-fit design--decision solution before the preferences defining that solution are even applied.

\noindent
This decision invalidity is demonstrated through a like-for-like comparison of Pareto-based selection and the Integrative Maximisation of Aggregated Preferences (IMAP) optimisation method across constructed multi-objective design--decision problems involving two system performances and up to four preference functions elicited from up to two actors. Both methods are applied to identical feasible solution spaces using identical actor-preference information. The problem set covers independent and interdependent relationships between feasible performance values, monotonic and non-monotonic preference functions, and conflicting actor preferences concerning the same performance. In a subset of these problems, Pareto filtering excludes the best-fit solution, whereas IMAP identifies it.

\noindent
This design--decision invalidity occurs even in these deliberately simplified problems. Real-world integrated design and group decision-making problems are substantially more complex, involving more system performances, non-monotonic preferences, performance interdependencies, and heterogeneous or conflicting actor interests. Pareto optimality is therefore insufficient to guarantee design--decision validity. Whereas Pareto analysis can support numerical ordering and exploration in performance space, IMAP establishes a preference-consistent design--decision process by modelling and aggregating actor preferences \emph{a priori} and maximising the resulting aggregate preference over the complete feasible and acceptable solution space.

\end{abstract}

\vspace*{1.0em}

\small
\noindent\textbf{Keywords:}
Pareto dominance concept;
Integrative Maximisation of Aggregated Preferences (IMAP);
Multi-objective design optimisation (MODO);
Multi-criteria decision-making (MCDM);
Design--decision systems;
Constraint Programming (CP);
Mixed-Integer (Non)Linear Programming (MI(N)LP);
Parametric Design;
Socio-technical Algorithmics;
Complex design--decision problems.
\\

\vspace{0.2em}

\normalsize
}

\section{Introduction}
\label{sec:intro}

Integrated design and group decision-making share a common task: identifying a feasible design alternative that best fits the interests of the actors involved. An alternative is a possible solution to the problem, described by system performances such as completion time and cost. These performances, also called objectives in many conventional optimisation approaches, allow alternatives to be compared but do not by themselves determine which alternative is best. That depends on the actors' preferences concerning different performance levels. The alternative that best fits their aggregated interests is referred to here as the \emph{best fit for common purpose}~\citep{Wolfert2026PreferenceBased,VanHeukelum2024}.

An actor's preference concerning a performance can be represented by a preference function that maps each feasible performance value to a degree of preference. In other words, a preference function represents how desirable an actor considers different levels of a performance. Integrating such functions into a design--decision problem requires Preference Function Modelling (PFM)~\citep{Barzilai2010Foundations}. Preference is subjective, relative, and context-dependent and is therefore measured on a one-dimensional affine scale that provides a basis for choice. Note that preference, utility, worth, and value are synonyms~\citep{Barzilai2022}. A preference function may be monotonically increasing or decreasing, but it need not be linear or even monotonic. It may represent acceptable or most-preferred ranges, thresholds, or saturation effects (e.g., system reliability or maintenance investments), where further improvements beyond an already highly preferred performance yield only marginal preference gains. Alternatively, it may be non-monotonic, with an interior maximum or minimum corresponding to a preferred or non-preferred target value (e.g., project completion time, urban density, staffing level, noise level, or operating temperature), where performance values on either side of the target level are respectively less preferred or more preferred than the performance of interest at that target level. Moreover, different actors may have conflicting preference functions concerning the same performance. In design--decision problems, preference conflicts can therefore arise not only between performances but also within a single performance~\citep{Wolfert2026PreferenceBased,TeuberWolfert2024ConfrontingConflicts}.

Thus, every multi-criteria design--decision method that ends in a choice relies on preference information, whether explicitly represented or implicitly assumed. If the applied preference information does not adequately represent the actors' interests, the resulting choice may be decision-invalid. Methods differ in where and how preference information enters the design--decision process~\citep{miettinen1999}. In Pareto-supported optimisation, a fixed improvement direction (minimisation or maximisation) is prescribed for each performance. Pareto filtering then orders the solution space by dominance to identify the non-dominated alternatives: those for which no feasible alternative improves at least one performance without worsening any other under the prescribed directions. These alternatives form a Pareto front~\citep{ehrgott2005,Martins2022}, depending on the min--max orientation. For example, with two performances, four orientations are possible: min/min, min/max, max/min, and max/max, specifying the prescribed improvement direction for each performance. As additional performances are introduced, the number of possible orientations grows exponentially. The alternatives on each front remain unranked. Identifying the best fit for common purpose therefore requires the designer or decision-maker to introduce additional actor-preference information in an \emph{a posteriori} selection step.

Such an \emph{a posteriori} selection step is often performed using a weighted-utility method, in which the performances are combined into a single function, for example through a weighted sum. However, this may produce non-unique and preference-inconsistent rankings~\citep{Wolfert2026UniquePreference}. In contrast, in preference-based \emph{a priori} methods, actor preferences are modelled and aggregated using PFM before optimisation. The Integrative Maximisation of Aggregated Preferences (IMAP) optimisation method, combined with PFM, directly maximises the aggregated preference over the complete feasible and acceptable design--decision space~\citep{Wolfert2026PreferenceBased,VanHeukelum2024,teuber2025odycon}. The resulting alternative is the \emph{preference optimum} and represents the preference-consistent best fit for common purpose. This distinction raises a fundamental question underlying Pareto-supported optimisation: whether an ordering constructed independently of the actors' preferences can identify their best fit for common purpose, even when preference information is applied \emph{a posteriori}.

\subsection{Pareto's Claim}
\label{sec:pareto-claim}

The view that alternatives can be ordered independently of the actors' preferences originates in Pareto's claim, as quoted in Barzilai \citet{Barzilai2010Foundations}, that ``the entire theory of economic equilibrium is independent of the notions of (economic) utility''. \citet{Barzilai2010Foundations} argues that this claim is incorrect: Pareto's development requires mathematical operations that are not supported by the ordinal utility scales on which the claim relies. \citet{Barzilai2022} discusses related errors in group decision theory.

Pareto dominance nevertheless remains central to multi-objective optimisation, including engineering design optimisation and evolutionary algorithms~\citep{Martins2022,Marler2004,Deb2002,Deb2014}. In \emph{a posteriori} approaches, the absence of actor-preference information during Pareto filtering is addressed by introducing it afterwards to select from the resulting Pareto front~\citep{miettinen1999,ehrgott2005}. Whether this subsequent use of preference information is sufficient to identify the preference optimum remains to be established and therefore constitutes the central research question of this paper.

\subsection{Further Limitations of Pareto Analysis}
\label{sec:relatedwork}

The limitations of Pareto analysis extend beyond Pareto's preference-independence claim. More broadly, \citet{french2023reflections} situates formal multi-criteria decision-making within an iterative process encompassing sense-making and modelling, analysis and exploration, and interpretation and implementation. From this overarching perspective, Pareto analysis supports the numerical ordering and exploration of alternatives in performance space but does not itself provide the preference measurement and aggregation required to identify a unique, actionable, and preference-consistent design--decision outcome.

More specifically, studies in the PFM tradition~\citep{Wolfert2026PreferenceBased,VanHeukelum2024,teuber2025odycon} argue that a Pareto front is insufficient as a basis for design--decision support. Pareto filtering identifies a subset of non-dominated alternatives based on prescribed performance-improvement directions but does not itself represent the actors' possibly non-linear, non-monotonic, or conflicting preference functions. Designers and decision-makers ultimately seek a design--decision configuration that best fits the actors' aggregated interests, rather than a subset of alternatives that are merely non-dominated in performance space. Preference-aware Pareto variants incorporate preferences into the selection of candidates but still aim to approximate a Pareto set~\citep{Ahmadianshalchi2023}. Additional preference information therefore remains necessary to move from numerical performance-space ordering to identifying the best fit for common purpose. This limitation is one of the principal motivations underlying the development of IMAP~\citep{Wolfert2026PreferenceBased,VanHeukelum2024,teuber2025odycon}.

A separate limitation identified in many-objective optimisation is that the proportion of non-dominated alternatives can increase with the number of objectives, thereby reducing the discriminative power of Pareto dominance~\citep{Purshouse2007}. The present study addresses a different unresolved question. Even if a preference-consistent ranking is applied \emph{a posteriori}, can it identify the preference optimum from the alternatives retained by Pareto filtering? Answering this question requires establishing whether the Pareto-filtered set necessarily contains the preference optimum.

\subsection{Aim and Contributions}
\label{sec:what-paper-shows}

This paper examines whether Pareto-based selection necessarily retains the preference optimum when actor preferences are non-monotonic, target-based, or conflicting. The central research question is: \textbf{Does the union of the Pareto fronts obtained from all possible minimisation--maximisation orientations necessarily contain the preference optimum when actor-preference information is applied \emph{a posteriori}?}

To address this question, Pareto-based selection and IMAP are compared across six constructed multi-objective design--decision problems. Each problem involves two system performances, up to two actors, and up to four preference functions. Together, the problems cover rectangular and circular solution spaces, uncoupled and coupled feasible performance values, monotonic and non-monotonic preference functions, and conflicting actor preferences concerning the same performance.

The comparison is performed on a like-for-like basis. All six problems involve two performances and at most two actors. For each problem, both methods use the same feasible and acceptable design--decision space, the same actor-preference information, and the same preference aggregation. The difference lies solely in when the resulting aggregated-preference function is applied. In Pareto-based selection, alternatives are first filtered by dominance under all four possible minimisation--maximisation orientations defined by the two performances, after which the aggregated-preference function is applied \emph{a posteriori}. In IMAP, the same aggregated-preference function is used \emph{a priori} as the optimisation criterion over the complete feasible and acceptable design--decision space.

By holding the solution space, actor-preference information, and preference aggregation constant, this construction isolates the effect of Pareto filtering. The six problems are not intended to estimate how frequently Pareto-based selection excludes the preference optimum in practice, but to determine whether it provides a general guarantee that the preference optimum remains available for subsequent preference-consistent selection. The paper thereby provides a direct comparison between preference-consistent \emph{a posteriori} selection from the union of all Pareto fronts and \emph{a priori} preference-based optimisation over the complete design--decision space.

The paper is structured as follows. Section~\ref{sec:comparison} introduces the common problem formulation, presents the six constructed problems, and compares the Pareto-based and IMAP outcomes. Section~\ref{sec:reflections} discusses the findings, their implications for real-world design--decision problems, the role of preference-guided search, and the resulting design--decision conclusions.

\section{Pareto-Based Selection versus IMAP for Identifying a Design--Decision Optimum}
\label{sec:comparison}

This section compares Pareto-based selection and IMAP across six constructed multi-objective design--decision problems. The six problems use canonical settings from multi-objective optimisation~\citep{ehrgott2005,Martins2022}, allowing Pareto dominance, Pareto fronts, and preference aggregation to be examined without unnecessary complexity and providing a transparent basis for comparing the two methods. Both methods are applied to the same feasible and acceptable design--decision space using the same actor-preference information and preference aggregation. Because IMAP maximises the aggregated preference over the complete space by construction, it identifies the preference optimum. The comparison therefore examines whether the numerically constructed Pareto fronts, using the same aggregated-preference function \emph{a posteriori}, also retain and identify this optimum. The preference optimum is labelled \emph{found} when it is retained on at least one of the constructed Pareto fronts and selected from their union, and \emph{missed} when it lies on none of these fronts and is therefore unavailable for subsequent selection.
The software implementations of the six problems are available from the repository specified in the Data Availability statement.

\subsection{Common Multi-Objective Design--Decision Problem Formulation}
\label{sec:setting}

The IMAP--Pareto comparison is conducted entirely in the two-dimensional performance space. For all six problems, the number of performances is $I=2$ and the number of actors is at most $K=2$. This restriction enables all $2^{I}=4$ Pareto orientations to be analysed exhaustively and the resulting preference landscapes, Pareto fronts, and preference optima to be visualised directly in two and three dimensions. The purpose of this simplified setting is analytical clarity rather than realism: it provides the minimal setting needed to examine whether Pareto-based selection guarantees retention of the preference optimum. A feasible and acceptable alternative is therefore represented directly by its performance vector $\mathbf{f}=(f_1,f_2)\in S$, where $S$ denotes the achievable set of feasible and acceptable performance vectors. The terms \emph{interior} and \emph{boundary} refer to this set. Each actor expresses a preference concerning each relevant performance through a preference function. The figures are computed from sampled performance vectors within each achievable set.

The problems use monotonic and non-monotonic preference functions. A monotonic preference function, whether linear or non-linear, represents a consistent preference for more or less of a performance and therefore defines one improvement direction over the feasible range. A non-monotonic preference function does not define one fixed improvement direction over this range. In particular, such a function may attain a global maximum at an interior performance value, with preference increasing towards this value from one direction and decreasing beyond it in the other.

\paragraph{Pareto criterion.} An orientation $\mathbf{e}\in\{-1,+1\}^{I}$ fixes a minimisation or maximisation direction for each performance, where $I$ is the number of performances, $e_i=+1$ denotes maximisation, and $e_i=-1$ denotes minimisation. Under orientation $\mathbf{e}$, a performance combination $\mathbf{a}\in S$ dominates another performance combination $\mathbf{b}\in S$ if

\begin{equation}
  \mathbf{a}\succ_{\mathbf{e}}\mathbf{b}
  \iff
  e_i a_i \geq e_i b_i \ \text{for all } i
  \quad\text{and}\quad
  e_j a_j > e_j b_j \ \text{for some } j,
  \label{eq:dominance}
\end{equation}

and the corresponding Pareto front is the set of performance combinations in $S$ that are not dominated under that orientation:

\begin{equation}
  \mathcal{P}_{\mathbf{e}}(S)
  =
  \left\{
    \mathbf{a}\in S:
    \nexists\,\mathbf{b}\in S
    \text{ such that }
    \mathbf{b}\succ_{\mathbf{e}}\mathbf{a}
  \right\}.
  \label{eq:front}
\end{equation}

In general, $I$ performances define $2^{I}$ possible orientations and corresponding Pareto fronts, although fronts obtained under different orientations need not be distinct. For the two performances considered here, the four orientations are min/min, min/max, max/min, and max/max. Pareto dominance and Pareto fronts are standard concepts~\citep{ehrgott2005,miettinen1999}.

\paragraph{IMAP criterion.} For $K$ actors and $I$ performances, let $z_{k,i}(f_i)$ denote the z-normalised preference of actor $k$ concerning performance $i$. Since these preference functions are represented on a common linear preference scale, they can be consistently aggregated using the linear weighted-centroid operation performed by the A-FINE aggregator $\mathbf{A}$~\citep{Wolfert2026PreferenceBased}. The software implementation is available through the repository specified in the Data Availability statement. The resulting aggregated-preference function is

\begin{equation}
  \label{eq:imap}
  Z(\mathbf{f})
  =
  \mathbf{A}\bigl(z_{k,i}(f_i);\,w'_{k,i}\bigr)
  =
  \sum_{i=1}^{I}A_i(\mathbf{f})
  =
  \sum_{i=1}^{I}\sum_{k=1}^{K}
  w'_{k,i}\,z_{k,i}(f_i).
\end{equation}

Here, $A_i(\mathbf{f})=\sum_{k=1}^{K}w'_{k,i}z_{k,i}(f_i)$ is the weighted aggregation of the relevant actors' preferences concerning performance $i$ and constitutes that performance's contribution to $Z(\mathbf{f})$. The effective contribution weight $w'_{k,i}=w_k w_{k,i}$ combines the global weight $w_k$ of actor $k$ with the actor's local weight $w_{k,i}$ for performance $i$.

IMAP models and aggregates the actor preferences \emph{a priori} and directly identifies the aggregated preference optimum 

\begin{equation}
  \label{eq:argmax}
  \mathbf{f}^{*} =
  \arg\max_{\mathbf{f}\in S} Z(\mathbf{f}).
\end{equation}

Thus, IMAP maximises the aggregated preference over the complete achievable set $S$, including all alternatives on the Pareto fronts.

\paragraph{Like-for-like comparison.}
Both methods use the same achievable performance set $S$, actor-preference functions, A-FINE aggregation, and resulting aggregated-preference function $Z(\mathbf{f})$. In IMAP, the preference functions are aggregated \emph{a priori}; for the like-for-like Pareto-based comparison, the same preference functions, weights, and A-FINE aggregation are applied \emph{a posteriori} to the alternatives retained by Pareto filtering. For this comparison, IMAP, like Pareto-based selection, is formulated directly in performance space as $P(\mathbf{f})$, rather than through its standard design--decision-space formulation using the composite mapping $P(\mathbf{f}(\mathbf{x}))$. IMAP maximises $Z(\mathbf{f})$ over the complete set $S$, including all alternatives on the Pareto fronts, whereas Pareto-based selection maximises $Z(\mathbf{f})$ only over the union of the four fronts. The preference optimum is \emph{found} if it lies on at least one front and is selected from their union; otherwise, it is \emph{missed}. Furthermore, three structural properties guide the like-for-like comparison:

\begin{itemize}[leftmargin=*]
  \item \emph{Separability.} Because $Z(\mathbf{f})$ is additively separable across performances, a change in $f_i$ affects $Z(\mathbf{f})$ only through $A_i(f_i)$. When only one actor expresses a preference concerning performance $i$, the shape of $A_i(f_i)$ is therefore determined by that actor's preference function.

  \item \emph{Direction.} Pareto dominance requires a fixed minimisation or maximisation direction for each performance. A monotonic preference function provides such a direction, whereas a non-monotonic preference function with an interior global maximum does not define one fixed direction over the complete feasible range. To avoid relying on an arbitrarily prescribed orientation, all $2^{I}$ possible orientations are therefore considered.

  \item \emph{Shape of $S$.} The geometry of the achievable set determines where the Pareto fronts lie. For the circular set considered here, the union of the four fronts contains the complete boundary. For the rectangular set considered here, each orientation yields one corner as its Pareto front; hence, the non-corner boundary points and all interior points lie on none of the four fronts.
\end{itemize}

\textsc{Note on problem numbering and terminology.} Each problem is identified by two numbers: the first denotes the problem number (1--6), and the second denotes the number of actors (1 or 2). Moreover, ``peaked at'' and ``has a target value at'' are used to mean ``attains a global maximum at''. The ``preference optimum'' is the alternative that attains the maximum aggregated preference and is therefore the best fit for common purpose.

\subsection{Uncoupled Rectangular Solution Space: Problems 1--2}
\label{sec:square}

For Problems 1--2, the solution space is the rectangle

\begin{equation}
  S_1
  =
  \left\{
    (f_1,f_2)\in\mathbb{R}^2:
    a\leq f_1\leq b,\;
    c\leq f_2\leq d
  \right\}
  =
  [a,b]\times[c,d].
  \label{eq:rectangle-solution}
\end{equation}

The performances are uncoupled: every feasible value of $f_1\in[a,b]$ can be combined with every feasible value of $f_2\in[c,d]$. The four Pareto fronts reduce to the corner points $(a,c)$, $(a,d)$, $(b,c)$, and $(b,d)$ for the min/min, min/max, max/min, and max/max orientations, respectively. Under each orientation, the corresponding corner dominates every other point in $S_1$.

\begin{figure}[H]
\centering
\begin{tikzpicture}[font=\small]
  \fill[black!7] (-1.7,-1.3) rectangle (1.7,1.3);
  \draw[black!45, line width=0.8pt] (-1.7,-1.3) rectangle (1.7,1.3);
  \draw[->] (-2.8,-2.2) -- (2.7,-2.2) node[below left] {$f_1$};
  \draw[->] (-2.8,-2.2) -- (-2.8,2.3) node[left] {$f_2$};
  \node[draw=black!80, fill=frontMaxMax, minimum size=6pt, inner sep=0pt] at (1.7,1.3) {};
  \node[draw=black!80, fill=frontMinMax, minimum size=6pt, inner sep=0pt] at (-1.7,1.3) {};
  \node[draw=black!80, fill=frontMinMin, minimum size=6pt, inner sep=0pt] at (-1.7,-1.3) {};
  \node[draw=black!80, fill=frontMaxMin, minimum size=6pt, inner sep=0pt] at (1.7,-1.3) {};
  \node[text=frontMaxMax, anchor=south] at (1.7,1.45) {max/max};
  \node[text=frontMinMax, anchor=south] at (-1.7,1.45) {min/max};
  \node[text=frontMinMin, anchor=north] at (-1.7,-1.45) {min/min};
  \node[text=frontMaxMin, anchor=north] at (1.7,-1.45) {max/min};
  \node[align=center, font=\footnotesize] at (0,0) {$\operatorname{int}S$ and edges:\\dominated under\\every orientation};
\end{tikzpicture}
\caption{Pareto ``fronts'' (corner points) under the four minimisation--maximisation orientations for the rectangular solution space $S_1$.}
\label{fig:rectangle-fronts}
\end{figure}
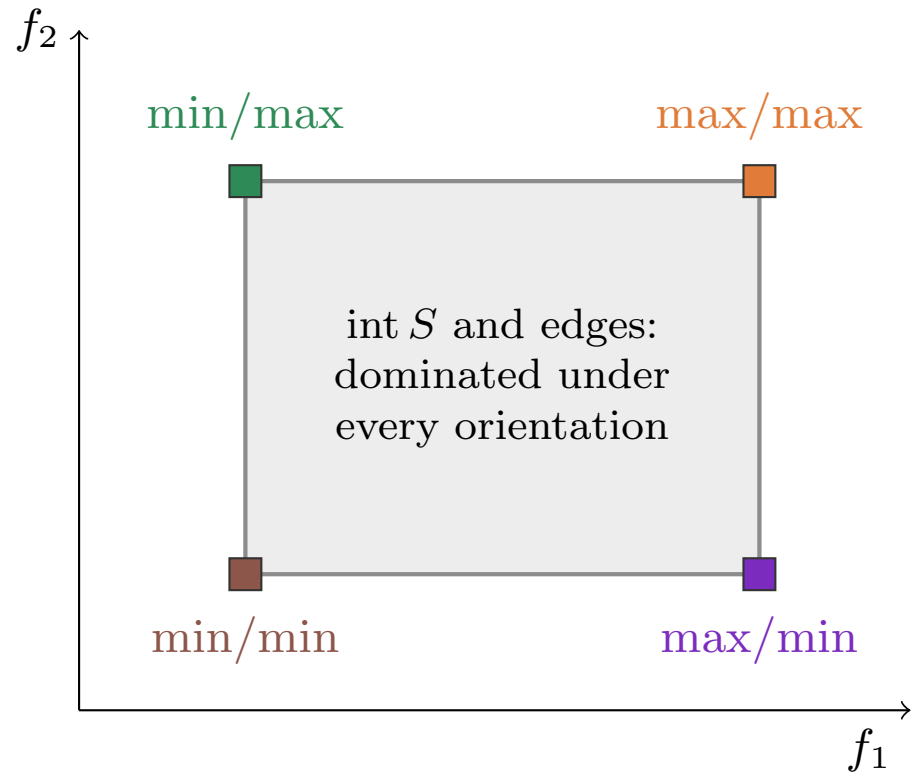

In Problems 1--2, one actor expresses preferences concerning both performances. With $K=1$ and $I=2$, the A-FINE aggregation in Eq.~\eqref{eq:imap} simplifies to

\begin{equation}
\label{eq:imap-one-actor}
Z(f_1,f_2)
=
A_1(f_1)+A_2(f_2)
=
w_1 z_1(f_1)+w_2 z_2(f_2).
\end{equation}

where $z_i(f_i)$ is the actor's z-normalised preference concerning performance $i$, $w_i$ is its effective contribution weight, and $A_i(f_i)=w_i z_i(f_i)$ is the corresponding contribution to $Z(f_1,f_2)$. The same formulation applies to the one-actor Problems 3--4.

\subsubsection{Problem 1-1: $P_1$, $P_2$ both linear}
Let the preference function for $f_1$ be linearly decreasing and that for $f_2$ linearly increasing, with equal local weights assigned to both performances. For one actor, each contribution $A_i(f_i)$ is a positive affine transformation of the corresponding preference function and is therefore strictly monotonic. Consequently, the preference optimum is $\mathbf{f}^{*}=(a,d)$, where $f_1$ attains its minimum and $f_2$ its maximum. This point coincides with the Pareto ``front'' of the min/max orientation, which follows directly from the two preference directions. Pareto-based selection therefore finds $\mathbf{f}^{*}$, and both methods agree. This result follows analytically from the separable linear formulation and is independent of sampling.

\begin{figure}[H]
\centering
\includegraphics[width=.75\textwidth]{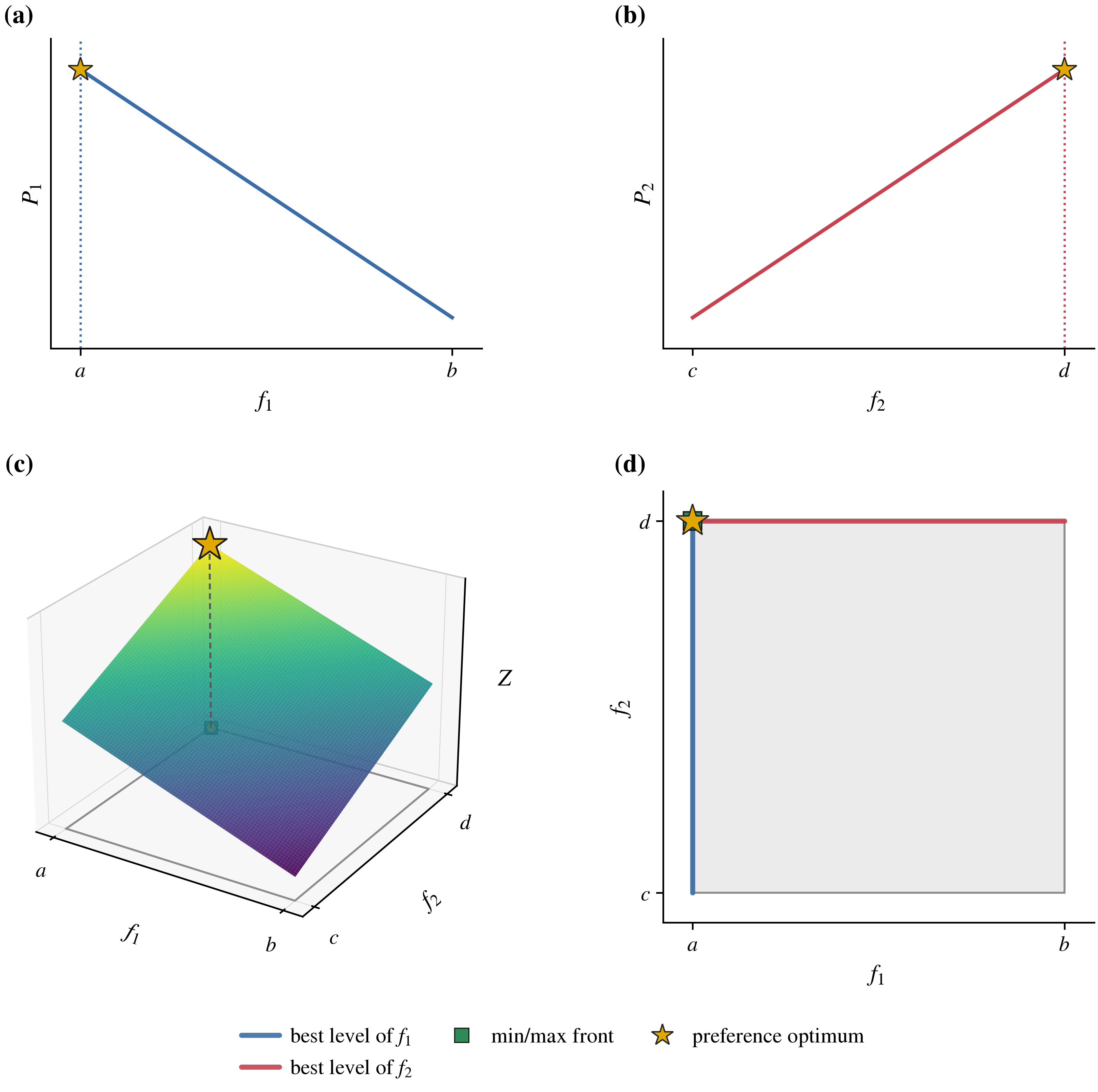}
\caption{Problem 1-1 with linear monotonic preferences and rectangular solution space $S_1$: (a,b) preference functions for $f_1$ and $f_2$; (c) aggregated-preference function $Z$, attaining its global maximum at the preference optimum $\mathbf{f}^{*}$; and (d) $\mathbf{f}^{*}$ coinciding with the min/max Pareto ``front'' (corner point).}
\label{fig:square-monotonic}
\end{figure}

\subsubsection{Problem 2-1: $P_1$ linear \& $P_2$ non-monotonic with interior global max}
\label{sec:square-peak}
The linearly decreasing preference function for $f_1$ is retained, while the preference function for $f_2$ is non-monotonic and attains its global maximum at an interior value $t$, where $c<t<d$. Consequently, $A_1$ is decreasing and $A_2$ attains its global maximum at $t$, yielding the preference optimum $\mathbf{f}^{*}=(a,t)$. This point lies on the left boundary of $S_1$, strictly between the corner points $(a,c)$ and $(a,d)$, and therefore lies on none of the four Pareto fronts. Because the preference function for $f_2$ does not define one fixed improvement direction over its complete feasible range, all four orientations are considered. Under each orientation, the corresponding corner dominates $\mathbf{f}^{*}$ but has a lower aggregated preference. Pareto-based selection therefore misses $\mathbf{f}^{*}$, whereas IMAP identifies it.

\begin{figure}[H]
\centering
\includegraphics[width=.75\textwidth]{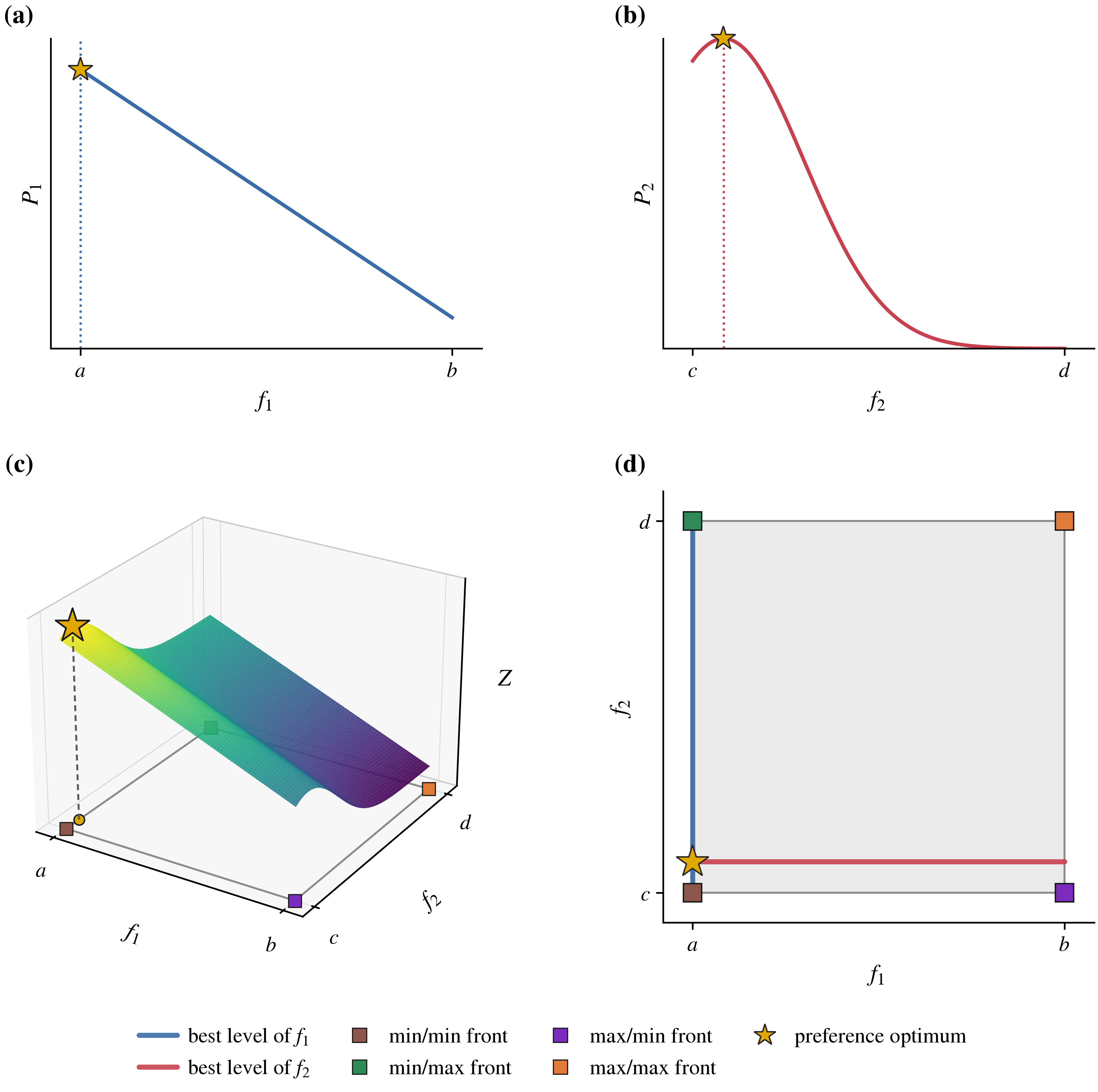}
\caption{Problem 2-1 with a linear decreasing preference for $f_1$, a non-monotonic preference for $f_2$ attaining an interior global maximum at $t$, and rectangular solution space $S_1$: (a,b) preference functions for $f_1$ and $f_2$; (c) aggregated-preference function $Z$ and preference optimum $\mathbf{f}^{*}$; and (d) the four Pareto ``fronts'' (corner points), the preferred performance levels, and $\mathbf{f}^{*}$ on the left boundary of $S_1$ between two corner points.}
\label{fig:one-case3}
\end{figure}

\subsection{Coupled Circular Solution Space: Problems 3--6}
\label{sec:circle}

For Problems 3--6, the achievable set is a circle inscribed in the box $[a,b]\times[c,d]$. The corresponding solution space is

\begin{equation}
  S_2=\Big\{\mathbf f:\ \Big(\frac{f_1-m_1}{h_1}\Big)^2+\Big(\frac{f_2-m_2}{h_2}\Big)^2\le 1\Big\},
  \label{eq:circle-solution}
\end{equation}

where $m_1=\tfrac{a+b}{2}$ and $m_2=\tfrac{c+d}{2}$ are the centre coordinates, and $h_1=\tfrac{b-a}{2}$ and $h_2=\tfrac{d-c}{2}$ are the corresponding half-widths of the box. For the circular solution space considered here, $h_1=h_2$. 

\begin{figure}[H]
\centering
\begin{tikzpicture}[font=\small]
  \fill[black!7] (0,0) circle[radius=1.6];
  \draw[->] (-2.8,-2.2) -- (2.7,-2.2) node[below left] {$f_1$};
  \draw[->] (-2.8,-2.2) -- (-2.8,2.3) node[left] {$f_2$};
  \draw[line width=1.6pt, frontMaxMax] (1.6,0) arc[start angle=0, end angle=90, radius=1.6];
  \draw[line width=1.6pt, frontMinMax] (0,1.6) arc[start angle=90, end angle=180, radius=1.6];
  \draw[line width=1.6pt, frontMinMin] (-1.6,0) arc[start angle=180, end angle=270, radius=1.6];
  \draw[line width=1.6pt, frontMaxMin] (0,-1.6) arc[start angle=270, end angle=360, radius=1.6];
  \node[text=frontMaxMax, anchor=south west] at (1.13,1.13) {max/max};
  \node[text=frontMinMax, anchor=south east] at (-1.13,1.13) {min/max};
  \node[text=frontMinMin, anchor=north east] at (-1.13,-1.13) {min/min};
  \node[text=frontMaxMin, anchor=north west] at (1.13,-1.13) {max/min};
  \node[align=center, font=\footnotesize] at (0,0) {$\operatorname{int}S$:\\dominated under\\every orientation};
\end{tikzpicture}
\caption{Pareto fronts under the four minimisation--maximisation orientations for the circular solution space $S_2$. Each front is a boundary arc; together, the four fronts form the complete boundary of $S_2$.}
\label{fig:circle-fronts}
\end{figure}
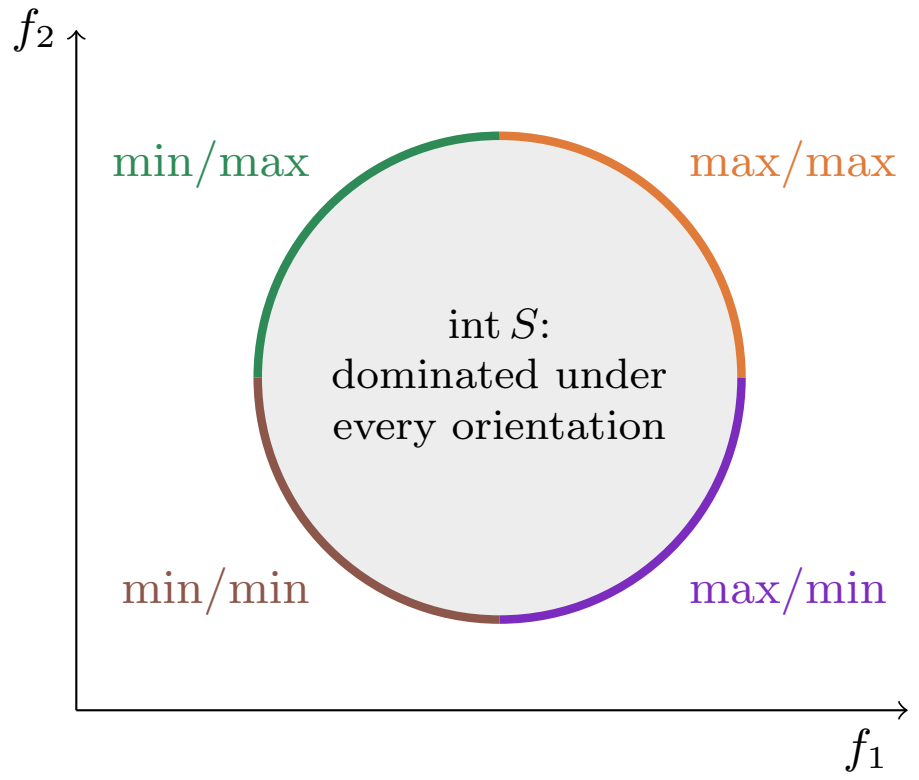

Thus, the bounding intervals have equal widths in the depicted performance space. This circular representation facilitates transparent geometric comparison and visualisation; the physical quantities, values, and units represented by $f_1$ and $f_2$ need not be identical.

The performances are coupled: the feasible values of $f_1$ and $f_2$ constrain one another, so not every combination within the bounding box is feasible. Each Pareto orientation defines one boundary arc. The min/min, min/max, max/min, and max/max fronts are the lower-left, upper-left, lower-right, and upper-right arcs, respectively. Together, the four fronts form the complete boundary of $S_2$, whereas the interior points lie on none of the four fronts. Problems 3--4 involve one actor and use the A-FINE aggregation in Eq.~\eqref{eq:imap-one-actor}. Problems 5--6 involve two actors; accordingly, the A-FINE aggregation becomes

\begin{equation}
\label{eq:imap-two-actors}
Z(f_1,f_2)
=
A_1(f_1)+A_2(f_2)
=
\sum_{k=1}^{2} w'_{k,1}z_{k,1}(f_1)
+
\sum_{k=1}^{2} w'_{k,2}z_{k,2}(f_2).
\end{equation}

\subsubsection{Problem 3-1: $P_1$ linear \& $P_2$ non-monotonic with interior global max}
\label{sec:circle-peak}

The linearly decreasing preference function for $f_1$ and the non-monotonic preference function for $f_2$ from Problem 2-1 are retained. Consequently, $A_1$ is decreasing, while $A_2$ attains its global maximum at the interior value $t$. Because the performances are coupled in $S_2$, the preference optimum $\mathbf{f}^{*}$ lies on the lower-left boundary arc, with $f_2$ close to but not equal to $t$. Moving along this arc towards a lower value of $f_1$ requires $f_2$ to move away from $t$, so $\mathbf{f}^{*}$ balances the resulting preference gain in $A_1$ against the preference loss in $A_2$. This boundary arc is the min/min Pareto front. Because the non-monotonic preference function for $f_2$ does not define one fixed improvement direction, all four orientations are considered. Their union contains the complete boundary of $S_2$ and therefore retains $\mathbf{f}^{*}$. Pareto-based selection consequently finds $\mathbf{f}^{*}$.

\begin{figure}[H]
\centering
\includegraphics[width=\textwidth]{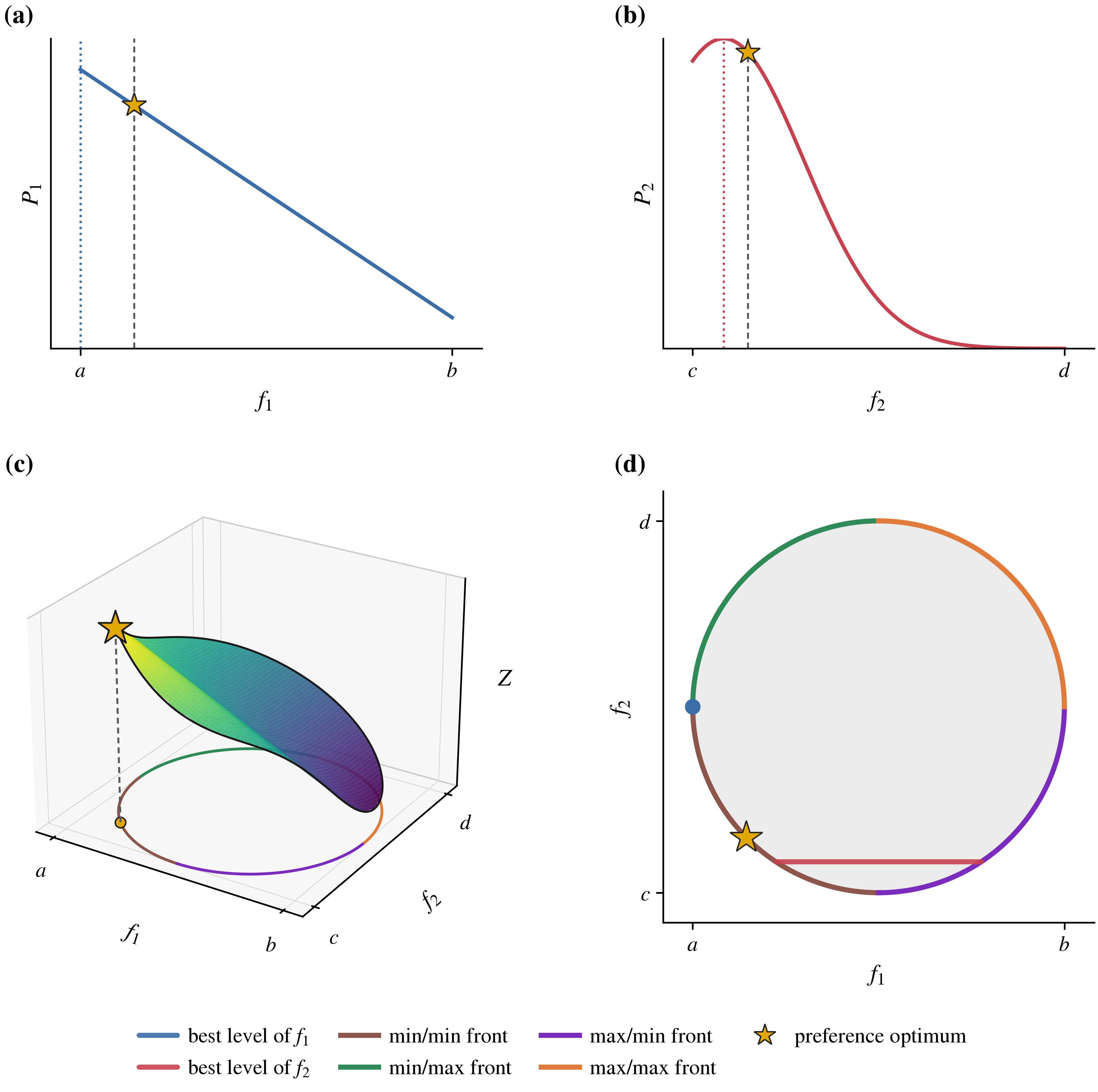}
\caption{Problem 3-1 with a linear decreasing preference for $f_1$, a non-monotonic preference for $f_2$ attaining an interior global maximum at $t$, and circular solution space $S_2$: (a,b) preference functions for $f_1$ and $f_2$; (c) aggregated-preference function $Z$, attaining its global maximum at the preference optimum $\mathbf{f}^{*}$; and (d) $\mathbf{f}^{*}$ located on the min/min Pareto front.}
\label{fig:one-case2}
\end{figure}

\subsubsection{Problem 4-1: $P_1$, $P_2$ both non-monotonic with interior global max}
\label{sec:circle-both}
Both preference functions are non-monotonic and attain their global maxima at interior values $t_1$ and $t_2$, with equal local weights assigned to the two performances. Consequently, $A_1$ and $A_2$ attain their global maxima at $t_1$ and $t_2$, respectively. Because $(t_1,t_2)\in\operatorname{int}(S_2)$, the preference optimum is $\mathbf{f}^{*}=(t_1,t_2)$. This interior point lies on none of the four Pareto fronts and is dominated under each orientation by at least one boundary point. However, every point in the union of the fronts has a lower aggregated preference than $\mathbf{f}^{*}$. Pareto-based selection therefore misses $\mathbf{f}^{*}$, whereas IMAP identifies it.

\begin{figure}[H]
\centering
\includegraphics[width=\textwidth]{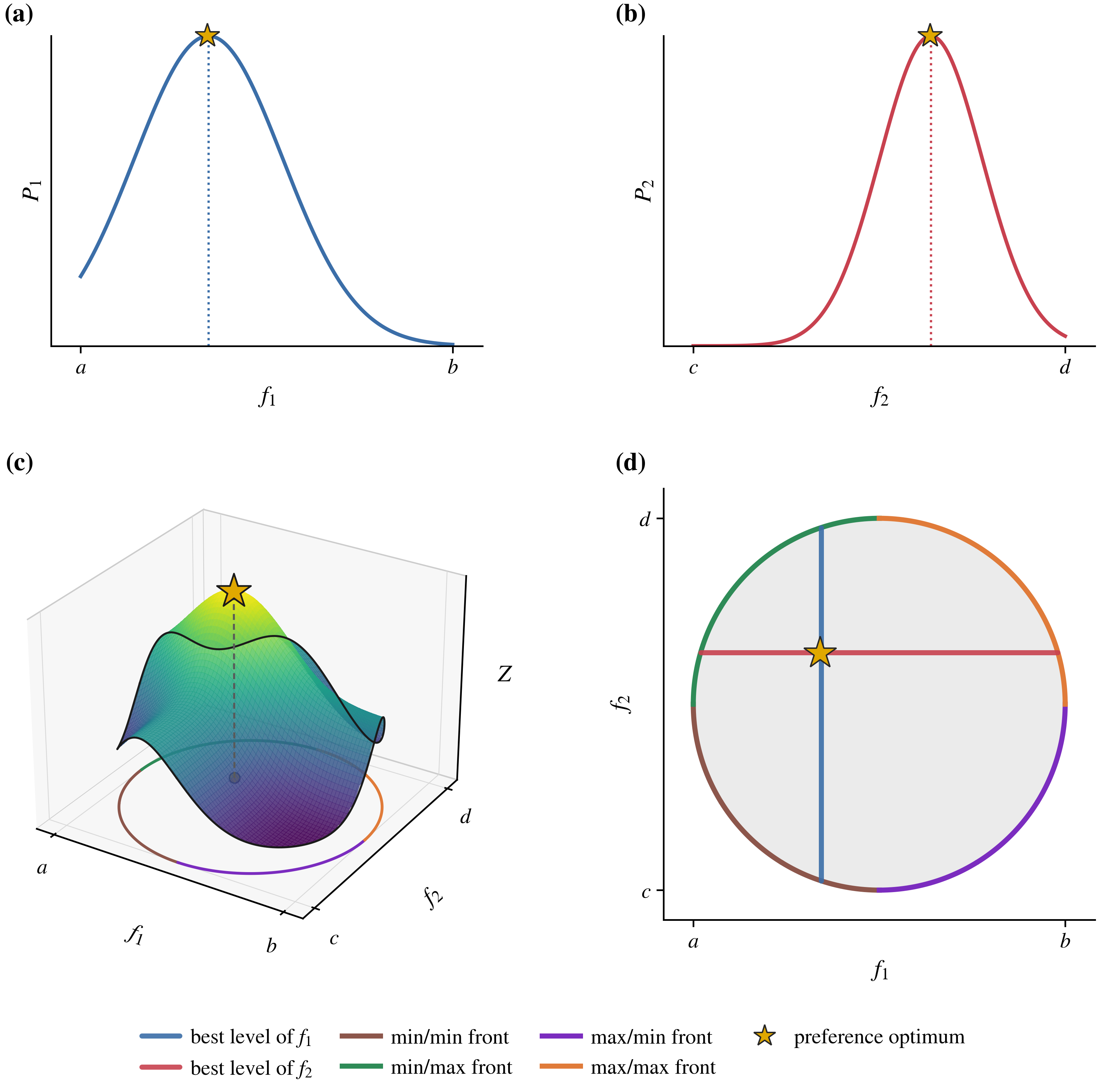}
\caption{Problem 4-1 with non-monotonic preferences for both $f_1$ and $f_2$, each attaining an interior global maximum, and circular solution space $S_2$: (a,b) preference functions for $f_1$ and $f_2$; (c) aggregated-preference function $Z$ and preference optimum $\mathbf{f}^{*}$; and (d) the four Pareto fronts, the preferred performance levels, and $\mathbf{f}^{*}$ in the interior of $S_2$.}
\label{fig:circle-both}
\end{figure}

\subsubsection{Problem 5-2: $P_{1,1}$, $P_{2,1}$ conflicting and linear \& $P_{1,2}$, $P_{2,2}$ equal non-monotonic with interior global max}
\label{sec:actors}

Two actors express preferences concerning both performances, with equal local weights assigned to $f_1$ and $f_2$. Their preferences concerning $f_2$ are identical and non-monotonic, each attaining its global maximum at the same interior value $t_2$. Consequently, the aggregated contribution $A_2$ also attains its global maximum at $t_2$. The actors differ only in their preferences concerning $f_1$: $P_{1,1}$ is linearly increasing, whereas $P_{2,1}$ is linearly decreasing.

Actor~1 has the larger global weight. The weighted aggregation of the opposing preferences therefore yields a linearly increasing contribution $A_1$. The global weights must be unequal for this problem: equal global weights would cause the opposing z-normalised linear preferences to cancel, making $A_1$ constant. The direction of $A_1$ is thus determined by the weighted aggregation of the actors' conflicting preferences, rather than by either actor's individual preference alone.

Because $A_1$ is increasing and $A_2$ attains an interior global maximum at $t_2$, the preference optimum $\mathbf{f}^{*}$ lies on the upper-right boundary arc of $S_2$, with $f_2$ close to but not equal to $t_2$. Moving along this arc towards a higher value of $f_1$ requires $f_2$ to move away from $t_2$, so $\mathbf{f}^{*}$ balances the preference gain in $A_1$ against the preference loss in $A_2$. This arc is the max/max Pareto front. Since the non-monotonic preferences concerning $f_2$ do not define one fixed improvement direction, all four orientations are considered. Their union retains $\mathbf{f}^{*}$, and Pareto-based selection therefore finds the same preference optimum as IMAP. If the actors' global weights were interchanged, $A_1$ would instead be decreasing and the preference optimum would lie on the upper-left, min/max Pareto front.

\begin{figure}[htbp]
\centering
\includegraphics[width=\textwidth]{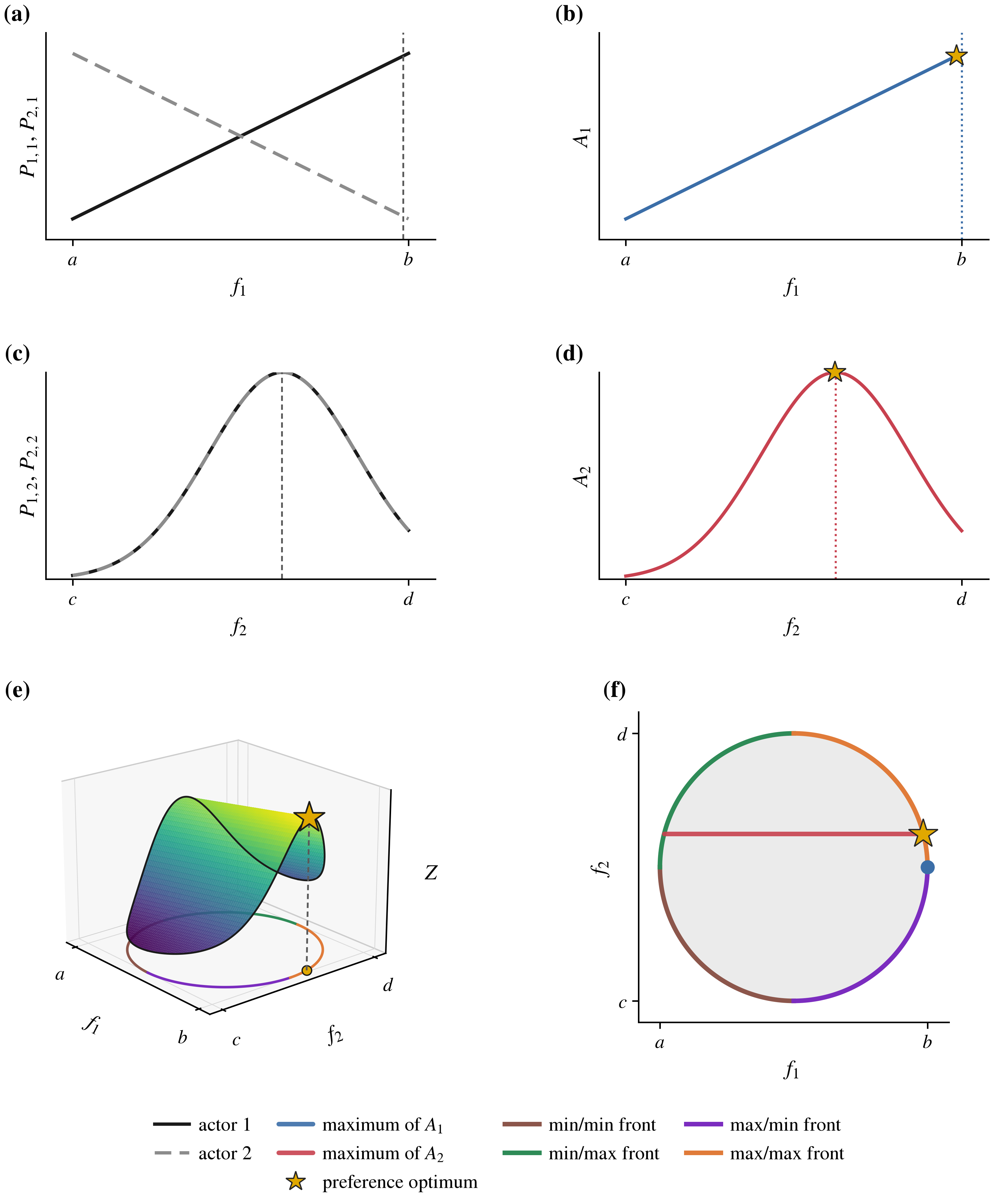}
\caption{Problem 5-2 with two actors having opposing linear preferences for $f_1$, identical non-monotonic preferences for $f_2$, and circular solution space $S_2$: (a,c) the actors' preference functions for $f_1$ and $f_2$; (b,d) the respective aggregated contributions $A_1$ and $A_2$, each attaining its global maximum at the indicated performance level; (e) aggregated-preference function $Z$ and preference optimum $\mathbf{f}^{*}$; and (f) the four Pareto fronts, the preferred performance levels, and $\mathbf{f}^{*}$ on the max/max front.}
\label{fig:two-case5}
\end{figure}

\subsubsection{Problem 6-2: $P_{1,1}$, $P_{2,1}$ conflicting and non-linear \& $P_{1,2}$, $P_{2,2}$ equal non-monotonic with interior global max}
\label{sec:case6}
Problems 5-2 and 6-2 form a controlled comparison: the actors' preference directions remain unchanged, while only the shapes of their preference functions concerning $f_1$ change from linear to non-linear. Actor~1's preference is increasing with diminishing returns, whereas Actor~2's preference is decreasing at an increasing rate. The actors' identical non-monotonic preferences concerning $f_2$ remain unchanged.

Although both preferences concerning $f_1$ are monotonic, their weighted aggregation $A_1$ is non-monotonic and attains its global maximum at an interior value $\hat{f}_1$. The aggregated contribution $A_2$ likewise attains its global maximum at the interior value $t_2$. Because $(\hat{f}_1,t_2)\in\operatorname{int}(S_2)$, the preference optimum is $\mathbf{f}^{*}=(\hat{f}_1,t_2)$.

The preference optimum lies on none of the four Pareto fronts and is dominated under each orientation by at least one boundary point. Nevertheless, every point in the union of the fronts has a lower aggregated preference than $\mathbf{f}^{*}$. Pareto-based selection therefore misses $\mathbf{f}^{*}$, whereas IMAP identifies it.

The comparison with Problem 5-2 isolates the effect of preference-function shape. Despite the unchanged individual preference directions, the aggregated contribution $A_1$ changes from monotonic to non-monotonic, and the preference optimum moves from the boundary to the interior of $S_2$. Consequently, the Pareto-based result changes from \emph{found} in Problem 5-2 to \emph{missed} in Problem 6-2.

\begin{figure}[htbp]
\centering
\includegraphics[width=\textwidth]{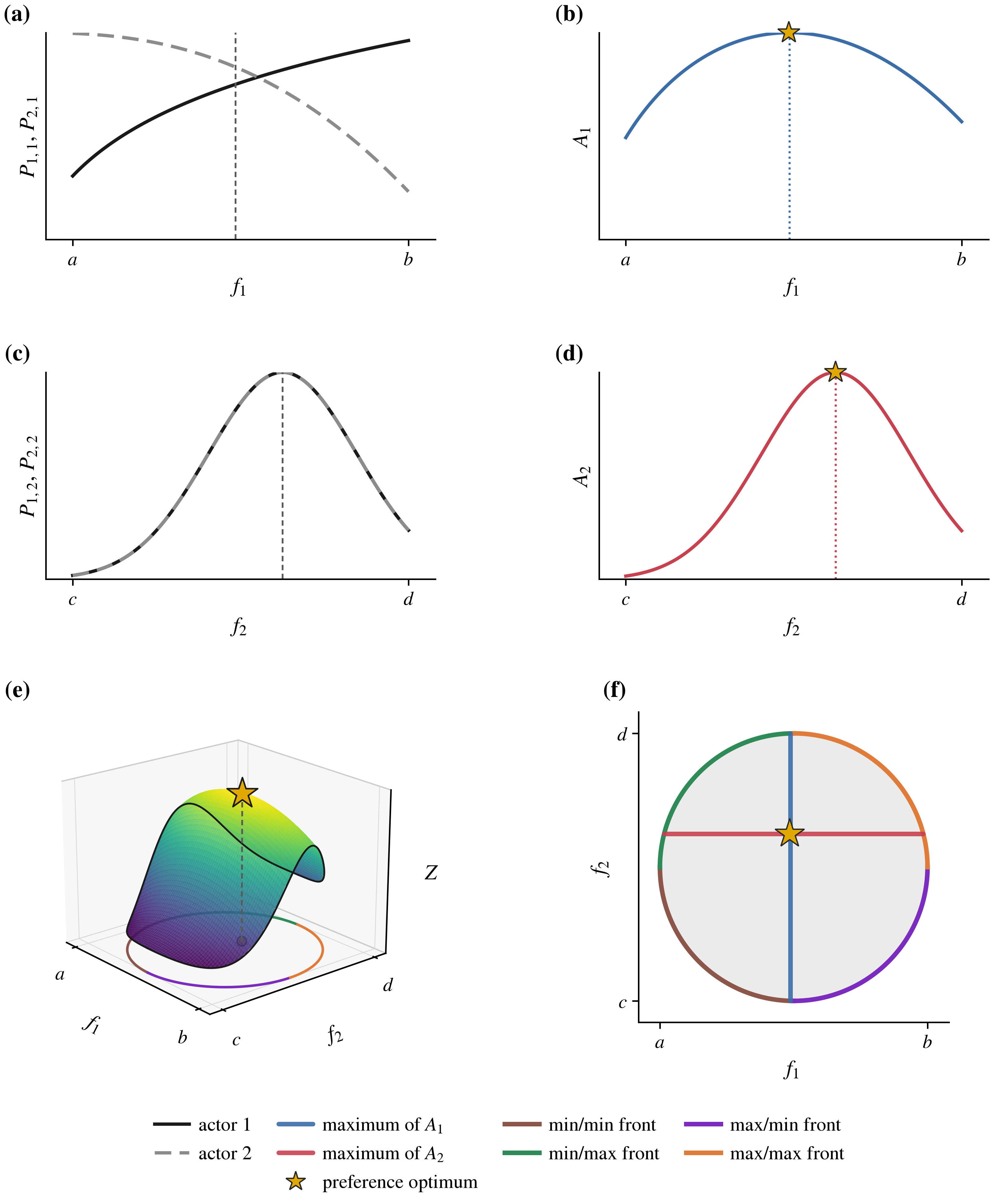}
\caption{Problem 6-2 with two actors having opposing non-linear preferences for $f_1$, identical non-monotonic preferences for $f_2$, and circular solution space $S_2$: (a,c) the actors' preference functions for $f_1$ and $f_2$; (b,d) the respective aggregated contributions $A_1$ and $A_2$, each attaining an interior global maximum; (e) the aggregated-preference function $Z$ and preference optimum $\mathbf{f}^{*}$; and (f) the four Pareto fronts, the most-preferred performance levels, and $\mathbf{f}^{*}$ in the interior of $S_2$.}
\label{fig:two-case6}
\end{figure}

\subsection{Comparative Overview of the Six Problems}
\label{sec:summary}

Table~\ref{tab:comparison} summarises the location of the preference optimum and the outcome of Pareto-based selection for the six problems. Pareto-based selection finds the preference optimum in Problems 1-1, 3-1, and 5-2, where the optimum lies on at least one Pareto front. It misses the preference optimum in Problems 2-1, 4-1, and 6-2, where the optimum lies on none of the four fronts.

\begin{table}[htbp]
\centering
\small
\setlength{\belowcaptionskip}{6pt}
\newcommand{\casesep}{\noalign{\vskip\aboverulesep\hbox to \linewidth{\color{black!45}\leaders\hbox to 3pt{\hss\vrule height 0.3pt width 1.5pt\hss}\hfill}\vskip\belowrulesep}}
\begin{tabularx}{\textwidth}{@{}>{\raggedright\arraybackslash}p{0.56\textwidth}>{\raggedright\arraybackslash}X>{\raggedright\arraybackslash}p{0.24\textwidth}@{}}
\toprule
Problem - \# of actors / types of \textit{P}-functions & IMAP & Pareto\\
\midrule
\multicolumn{3}{@{}l}{\textit{$S_1$ - Rectangular}} \\
\addlinespace[2pt]
1-1 / $P_1$, $P_2$ both linear & corner point & found (min/max front) \\
\casesep
2-1 / $P_1$ linear \& $P_2$ non-monotonic (interior max) & edge point & missed \\
\midrule
\multicolumn{3}{@{}l}{\textit{$S_2$ - Circular}} \\
\addlinespace[2pt]
3-1 / $P_1$ linear \& $P_2$ non-monotonic (interior max) & boundary point & found (min/min front) \\
\casesep
4-1 / $P_1$, $P_2$ both non-monotonic (interior max) & interior point & missed \\
\casesep
5-2 / $P_{1,1}$, $P_{2,1}$ conflicting and linear \& \\$P_{1,2}$, $P_{2,2}$ equal non-monotonic (interior max)  & boundary point & found (max/max front) \\
\casesep
6-2 / $P_{1,1}$, $P_{2,1}$ conflicting and non-linear \& \\$P_{1,2}$, $P_{2,2}$ equal non-monotonic (interior max)  & interior point & missed \\
\bottomrule
\end{tabularx}
\caption{Results for the six design--decision problems: the location of the IMAP preference optimum $\mathbf{f}^{*}$ in $S_i$ and whether Pareto-based selection recovers it as the design with the highest aggregated preference $Z$. The relevant front is shown in parentheses when $\mathbf{f}^{*}$ is found; ``missed'' indicates that $\mathbf{f}^{*}$ lies on none of the four possible Pareto fronts.}
\label{tab:comparison}
\end{table}

In the rectangular solution space $S_1$, the corner optimum of Problem 1-1 coincides with the min/max front and is found, whereas the non-corner boundary optimum of Problem 2-1 lies on none of the four corner fronts and is missed. In the circular solution space $S_2$, the boundary optima of Problems 3-1 and 5-2 lie on the min/min and max/max fronts, respectively, and are found. The interior optima of Problems 4-1 and 6-2 lie on none of the four fronts and are missed. Problems 5-2 and 6-2 provide the controlled comparison. The individual preference directions are unchanged, but replacing the conflicting linear preferences concerning $f_1$ with conflicting non-linear preferences changes the aggregated contribution $A_1$ from monotonic to non-monotonic. Correspondingly, the preference optimum moves from the boundary in Problem 5-2 to the interior in Problem 6-2, and the Pareto-based outcome changes from \emph{found} to \emph{missed}.

Across these six simplified, canonical problems, Pareto-based selection misses the preference optimum when at least one actor has a non-monotonic preference function with an interior global maximum (Problems 2-1 and 4-1), or even when conflicting monotonic actor preferences concerning one performance aggregate into a non-monotonic preference function while the aggregated preference concerning the other performance is also non-monotonic (Problem 6-2). These findings demonstrate that Pareto-based selection does not, in general, guarantee a decision-valid multi-criteria outcome.

\section{Conspection}
\label{sec:reflections}

The findings from the six constructed problems must be considered in the broader context of real-world design--decision optimisation. The failure of Pareto-based selection in three of these deliberately simplified, canonical problems does not provide an estimate of its failure rate in practice, but it disproves any general guarantee that Pareto filtering retains the preference optimum. If this guarantee already fails in a two-dimensional setting with at most two actors, four Pareto orientations, and explicitly defined preference functions, the substantially greater complexity of real-world design--decision problems provides no basis for assuming that it will be restored. The concluding subsection returns to the central research question.  

\subsection{Real-World Conditions Beyond the Constructed Problems}
\label{sec:reflection}

Real-world design--decision problems are substantially more complex than the constructed problems. Practical preference functions are seldom exclusively linear or monotonic; they may contain acceptable or most-preferred ranges, thresholds (e.g., budget caps and minimum availability requirements), saturation effects (e.g., system reliability and maintenance investments), and interior global maxima or minima (e.g., project completion time, urban density, staffing level, noise level, or operating temperature). Collectively, these characteristics demonstrate that the monotonic ``less-is-better'' and ``more-is-better'' preference representations implicitly assumed in Pareto analysis are often insufficient for representing real-world actor preferences. Multiple actors may also express heterogeneous or conflicting preferences concerning the same performance~\citep{Wolfert2026PreferenceBased,TeuberWolfert2024ConfrontingConflicts}. Accordingly, one fixed minimisation or maximisation direction cannot generally represent the complete preference structure associated with a performance in design--decision optimisation.

In practical applications, performances are often functions of shared design--decision variables, $\mathbf{f}(\mathbf{x})$, and are therefore commonly interdependent. Changing one design--decision variable may affect several performances simultaneously, producing coupled, non-linear, and constrained feasible spaces. The actors' performance-to-preference functions must consequently be evaluated through the composite mapping $P(\mathbf{f}(\mathbf{x}))$ and aggregated across performances and actors. The resulting best-fit design--decision point is obtained by maximising the aggregated preference over the feasible and acceptable design--decision space. In general, this optimum need not coincide with the maximum of any individual preference function, because it emerges from maximising the aggregation of the preferences of multiple actors across multiple interdependent performances within the coupled design--decision space. Pareto dominance does not represent this performance-to-preference and preference-aggregation structure; it orders feasible performance combinations numerically under prescribed improvement directions and therefore cannot guarantee decision-valid outcomes.

Increasing dimensionality introduces two further limitations. For $I$ performances, the number of possible minimisation--maximisation orientations grows exponentially as $2^{I}$, making exhaustive construction of all corresponding Pareto fronts rapidly impractical. Simultaneously, the proportion of mutually non-dominated alternatives generally increases with the number of performances, reducing the discriminative power of Pareto dominance~\citep{Purshouse2007,Deb2014}. Pareto analysis may therefore require an impractical number of orientations while producing increasingly large sets of non-dominated alternatives. Neither exhaustive orientation analysis nor enlarged non-dominated sets provide the preference measurement and aggregation required to identify the best fit for common purpose.

These real-world conditions provide no basis for expecting the retention guarantee absent from the simplified problems to reappear. Rather, non-linear and conflicting preferences, coupled performances, complex feasible spaces, increasing numbers of orientations, and declining dominance discrimination further separate numerical ordering in performance space from preference-consistent choice in design--decision space.

\subsection{Preference-Guided IMAP Design--Decision Optimisation}
\label{sec:choosing-imap}

Preference-guided IMAP design--decision optimisation reverses the conventional Pareto-based sequence. Pareto-based selection first filters the achievable set using numerical dominance and only then applies preference information to the retained alternatives. It can therefore identify the highest-preference alternative within the Pareto-filtered set but cannot guarantee that this alternative is the preference optimum over the complete feasible and acceptable design--decision space.

IMAP avoids this Pareto-filtering limitation by modelling and aggregating actor preferences before optimisation and using the resulting aggregated-preference function as the optimisation criterion from the outset. This requires no more preference information than a complete \emph{a posteriori} procedure, which ultimately also requires such information to select from the Pareto set. The decisive difference is that IMAP applies preference information before alternatives are discarded, rather than afterwards, when only the retained alternatives remain available for selection.

IMAP represents actor preferences through performance-to-preference functions, aggregates them in a common preference domain using the A-FINE aggregator, and maximises the resulting aggregated preference over the complete feasible and acceptable design--decision space~\citep{Timp2026,VanHeukelum2024,Wolfert2026PreferenceBased,teuber2025odycon}. Through the composite mapping from design--decision variables to performances and from performances to preferences, IMAP accommodates coupled performances, multiple actors, non-linear and non-monotonic preference functions, and conflicting actor preferences concerning the same performance.

IMAP is therefore recommended for preference-consistent design--decision optimisation. Numerical optimisation, approximation, or surrogate modelling may still support the search, but aggregated actor preference should define the optimisation criterion from the outset. Pareto filtering should not precede this preference-guided search, because it cannot guarantee retention of the preference optimum.

\subsection{Conclusions}
\label{sec:conclusions}

Pareto's preference-independence claim that alternatives can be ordered independently of preference information has already been shown to be incorrect. This study establishes the stronger result that subsequently applying the same actor-preference information \emph{a posteriori} does not restore a general guarantee of decision validity, because the preference optimum may already have been excluded by the preceding numerical dominance filter. This remains the case even when the union of the Pareto fronts under all possible minimisation--maximisation orientations is considered.

As a preference-independent method, Pareto analysis is therefore limited to numerical ordering and exploration in performance space. Pareto non-dominance does not imply preference optimality, and preference-consistent selection from a Pareto set does not guarantee identification of the best fit for common purpose. Consequently, Pareto-based selection cannot, in general, provide trustworthy, decision-valid multi-criteria decision support.

Preference-consistent design--decision optimisation requires actor preferences to define the optimisation criterion from the outset. IMAP implements this principle by modelling and aggregating actor preferences \emph{a priori} and maximising the resulting aggregated preference over the complete feasible and acceptable design--decision space. The overarching conclusion is therefore unequivocal: numerical Pareto ordering cannot replace preference-consistent deciding, nor can its design--decision deficiency be corrected by introducing preference information only after Pareto filtering. As a final caution, decision-makers should not entrust conventional Pareto-based approaches, including contemporary LLM-supported Pareto implementations, with identifying the best fit for common purpose. Numerical non-dominance is not a trustworthy substitute for preference-consistent choice and may exclude the best fit for common purpose, potentially leading to suboptimal design decisions.

\section*{Data Availability}

The current Preferendus software environment includes Preferendus-IMAP-IGS, which uses the A-FINE aggregator's decision-valid preference-aggregation method to support preference-consistent design--decision optimisation. The software is publicly available through the Preferendus GitHub organisation:

\url{https://github.com/Preferendus}.

The implementations of Problems 1--6 presented in this paper are available from the following GitHub repository:

\url{https://github.com/Preferendus/IMAP-vs-Pareto}.

\end{document}